\documentclass[a4paper,12pt]{article}

\title{Lower bounds on the Calabi functional}
\author{S. K. Donaldson}

\begin{document}
\maketitle
\newcommand{\bC}{\mbox{${\bf C}$}}
\newcommand{\bP}{\mbox{${\bf P}$}}
\newcommand{\uA}{\mbox{${\underline{A}}$}}
\newcommand{\uM}{\mbox{${\underline{M}}$}}
\newcommand{\um}{\mbox{${\underline{m}}$}}
\newtheorem{thm}{Theorem}
\newtheorem{prop}{Proposition}
\newtheorem{lem}{Lemma}

\section{Introduction}
A cornerstone of Atiyah and Bott's treatment \cite{kn:AB1} of Yang-Mills
theory over Riemann surfaces is a relation they discovered between the Yang-Mills functional
and  filtrations of a holomorphic bundle.  The relation can be stated
 as follows. Let $\Sigma$
be a compact Riemann surface with a fixed compatible metric, normalised
to have area $4\pi^{2}$. Let $E$ be a holomorphic vector bundle over $\Sigma$
and consider a flag ${\cal F}$ of sub-bundles
$$         0=E_{0}\subset E_{1}\dots \subset E_{q}=E .$$

Recall that for any bundle $V$ over $\Sigma$ the {\it slope} $\mu(V)$ of
 $V$ is defined to be the quotient of the degree of V by the rank.
 We say  the flag ${\cal F}$ is {\it slope-decreasing} if $\mu(E_{1})>
  \mu(E_{2})\dots>
 \mu(E)$. Equivalently the quotients $Q_{i}=E_{i}/E_{i-1}$ have
 $\mu(Q_{1})>\dots > \mu(Q_{q})$.  Define
 \begin{equation}\label{eq:phidef}
   \Phi({\cal F}) =\left(\sum_{i=1}^{q} \mu(Q_{i})^{2} {\rm rank}(Q_{i})\right)^{1/2}. \end{equation}
   Atiyah and Bott relate this algebro-geometric quantity to the Yang-Mills
 functional---the $L^{2}$ norm of the curvature $F_{A}$---restricted to
 compatible unitary connections $A$ on $E$. Their result is
 \begin{equation}\label{eq:ABrel}   {\inf}_{A} \Vert F_{A}\Vert = 
 {\sup}_{{\cal F},{\rm dec}} \Phi({\cal F})
 \end{equation} 
 where on the right hand side the supremum runs over the slope-decreasing
 flags. In fact this supremum is attained by the canonical
 Harder-Narasimhan filtration of $E$. The infimum on the left hand side
 is not in general attained: this happens if and only if $E$ is a direct
 sum of stable bundles. Notice that there is an easy  lower bound
 $ \Vert F_{A}\Vert\geq \mu(E) \sqrt{ {\rm rank}(E)}$ deriving from the
 fact that  $\frac{i}{2\pi} {\rm Tr}(F_{A})$ represents $c_{1}(E)$, This
 trivial lower bound is just $\Phi({\cal F}^{0})$ where  ${\cal F}^{0}$
 is the trivial flag (with
 $q=1$).
  So
 ${\sup}_{{\cal F},{\rm dec}} \Phi({\cal F}) \geq \mu(E) \sqrt{{\rm rank}\
 E}$, and it is easy to see that equality holds if and only if 
the bundle $E$ is semi-stable, i.e if there is no nontrivial slope-decreasing
flag.

The Atiyah-Bott result can be viewed as two statements: $\inf_{A} \geq
\sup_{{\cal F}}$ and $\inf_{A}\leq \sup_{{\cal F}}$. The proof of the first
of these involves a simple differential-geometric argument, turning on
the principle that \lq\lq curvature decreases in holomorphic sub-bundles
and increases in holomorphic quotients''(\cite{kn:AB1}, Proposition 8.13
and the remark in the second paragraph of page 575). The proof of the second relies
on
the theorem of Narasamihan and Seshadri on the existence of projectively
flat unitary connections on stable bundles.

In \cite{kn:Cal}, Calabi began the study of the $L^{2}$ norm of the scalar
curvature of Kahler metrics, running over a fixed Kahler class on a compact
Kahler manifold. This functional is equivalent to the $L^{2}$-norm of the
full curvature tensor, in that the two differ by topological terms. The
purpose of this paper is to establish an analogue of (one half of) the
Atiyah-Bott result for this Calabi functional (and some variants for
 $L^{p}$ norms).  Our result bears on the
algebraic case, so we suppose that $X$ is a smooth complex projective variety,
that $L\rightarrow X$ is a fixed ample line bundle and we consider Kahler
metrics $\omega$ in the class $c_{1}(L)$. We use the notion, introduced
in \cite{kn:D2}, of a {\it test configuration} ${\cal X}$ for
$X$. The detailed definition will be reviewed below, but in essence this
  comprises a ${\bf C}^{*}$-equivariant family $\pi:{\cal X}\rightarrow
{\bf C}$ with generic fibre $ X_{t}=\pi^{-1}(t)$ isomorphic to $X$, for $t\neq
0$. The central fibre $X_{0}$ need not be isomorphic to $X$, it may be a highly
singular variety or even a scheme,  but it has a ${\bf
C}^{*}$-action and the definition requires a lift of this action to a line bundle
${\cal L}$. We will define a numerical invariant $\Psi({\cal X})$ using the ${\bf
C}^{*}$ action on the vector spaces $H^{0}(X_{0}, {\cal L}^{k})$, related
to the generalised Futaki invariant. In stating our results, it is
convenient to work with a quantity $S(\omega)$ defined to be
$\frac{1}{4\pi}$ times the usual scalar curvature of the metric $\omega$.
(Thus $ S(\omega)\frac{\omega^{n}}{n!}= \frac{1}{2} 
\rho\wedge \frac{\omega^{n-1}}{(n-1)!}$, where $\rho$ is the Ricci form
representing $c_{1}(X)$.)
 Then our result
takes the form
\begin{equation}\label{eq:res}   \inf_{\omega} \Vert S(\omega)\Vert_{L^{2}} \geq \sup_{{\cal X}} \Psi({\cal
X}). \end{equation} 

There is a rather trivial lower bound on the Calabi functional, deriving
from the fact that the integral of $S(\omega)$ yields 
$ \langle \frac{1}{2(n-1)!} c_{1}(X) \cup \omega^{n-1}, [X]\rangle, $
so the average value $\hat{S}$ is a topological invariant of
the Kahler class. Thus
\begin{equation}\label{eq:lowbound} 
 \Vert S(\omega)\Vert_{L^{2}}^{2} = \Vert S-\hat{S}\Vert_{L^{2}}^{2} +
 \hat{S}^{2} \ {\rm Vol}(X) \geq \hat{S}^{2} \ {\rm Vol}(X). \end{equation}
 The essential feature of our definition of $\Psi({\cal X})$ is that we
get a better bound than this precisely when $X$ is not \lq\lq 
K-semistable'' in the sense of \cite{kn:D2}, just as in the 
Atiyah-Bott case. (There is also an extension of this discussion 
to extremal
metrics, see Section 2 below.) Thus one immediate
consequence of our Theorem (as Richard Thomas pointed out to the author)
 is a new and simpler proof of the fact that
the existence of  a constant scalar curvature Kahler metric
 implies K-semistability. This argument does not need any
  special constructions (as in \cite{kn:M}) in the case when the manifold has holomorphic vector
  fields. 
  
 It is natural to conjecture
that in fact
$$   \inf_{\omega} \Vert S(\omega)\Vert_{L^2} = \sup_{{\cal X}} \Psi({\cal
X}). $$    
This would be a variant of the conjectures in \cite{kn:Y},\cite{kn:T2},
\cite{kn:D2} relating K-stability
to the existence of constant scalar curvature metrics (and the extension
of that discussion to extremal metrics in \cite{kn:Sz}). However neither
conjecture immediately implies the other. 

We will explain the analogy with the Atiyah-Bott result in somewhat more
detail in Section 2. Both set-ups fit, formally, into the more general framework of
\lq\lq moment maps and stability'' discussed in \cite{kn:D1} for example.
The results can be seen as infinite-dimensional versions of part of the
theory developed by 
 Kirwan \cite{kn:K} in finite dimensions. The method of proof we use here
 is that of finite-dimensional approximation, in the mould of \cite{kn:D3},
 \cite{kn:D4}. That is, we derive our inequalities by studying
 the asymptotics of finite-dimensional problems, essentially of the kind
 considered by Kirwan. However the proof of the result in this paper
  is substantially simpler than that in \cite{kn:D3}. In particular we do not
 need to make use of the notion of \lq\lq balanced'' metrics.
 The essential ingredient in our proof is the asymptotic expansion for
 the \lq\lq density of states'' function due to Yau, Tian, Zelditch,Liu, Catlin and
 Ruan \cite{kn:T1}, \cite{kn:Ze},\cite{kn:L}, \cite{kn:C}, \cite{kn:Ru}. We also need some discussion of singular varieties and schemes 
  to allow us to apply the general moment map theory. The author first
  saw that this could be done using the point of view developed in the
   work of Zhang \cite{kn:Z}, Phong and Sturm
  \cite{kn:PS1},\cite{kn:PS2} involving the \lq\lq Chow norm'' and  an action on Chow vectors.
  However it turns out that one can avoid appealing to these concepts explicitly
  and we only need some comparatively straightforward technical facts to handle these issues
  of singularities (notably a result related to the 
  equivariant Riemann-Roch Theorem,
  Proposition 3 below).

 In 1997, and intermittently since, the author has discussed with
 X-X. Chen the possibility of obtaining lower bounds on the Calabi functional
 using a similar pattern of argument, but replacing the
  finite-dimensional approximations by the use of \lq\lq geodesic rays''
  in the space of Kahler potentials, in the manner of \cite{kn:Chen}, \cite{kn:CT}.
 This would have independent interest and Chen and the author hope
  to discuss this work in a a future article.
 
 The author is grateful to Richard Thomas and Xiu-Xiong Chen for discussions
 of these topics.
 
\section{Test configurations}
We begin by recalling the definition of a test configuration from \cite{kn:D1}.
Given an ample line bundle $L\rightarrow X$ a test configuration for the
pair $(X,L)$ consists of:
\begin{itemize}
\item a scheme ${\cal X}$ with a $\bC^{*}$ action;
\item  a flat
$\bC^{*}$-equivariant map $\pi:{\cal X}\rightarrow {\bC}$, with fibres
$X_{t}$;
\item an equivariant line bundle ${\cal L}\rightarrow {\cal X}$, ample
on all fibres;
\item for some $r>0$, an isomorphism of the pair $(X_{1}, {\cal L}\vert_{X_{1}})$ with the
original pair $(X,L^{r})$.
\end{itemize}
    
Thus we have a sequence of vector spaces $U_{k}=H^{0}(X_{0}, {\cal L}\vert_{X_{0}}^{k})$ with
$\bC^{*}$-actions. 
 Let
$  A_{k}: U_{k}\rightarrow U_{k}$ 
be the endomorphisms generating these actions (so $e^{t}\in {\bf C}^{*}$
acts as $e^{tA_{k}}$ on $U_{k}$). 
We are interested in the asymptotic behaviour, with respect to $k$,
 of the dimension of $U_{k}$
and the total weight of the action, i.e. the trace of the $A_{k}$.
These are given by Hilbert polynomials:
\begin{equation}\label{eq:aser} {\rm dim}\ U_{k}= a_{0} k^{n} + a_{1} k^{n-1} + \dots, \end{equation}

\begin{equation}\label{eq:bser}  {\rm Tr} A_{k}= b_{0} k^{n+1} + b_{1} k^{n} + \dots, \end{equation}
for large $k$, as discussed in \cite{kn:D1}, \cite{kn:RT}. We go one step
further and define a positive number $Q$ by the leading term in the polynomial
defining the trace of the squares:\begin{equation}\label{Qdef}
   {\rm Tr} A_{k}^{2} \sim Q\  k^{n+2}.  \end{equation}
 It is clear that that $Q\geq 0$, one way of seeing that $Q$ is strictly
 positive is to use (20) in Section 5 below.
 Notice that, from general theory, the $a_{i}, b_{i}$ and $Q$ are rational
 numbers. 

Now we define
\begin{equation}\label{eq:psidef}     
 \Psi({\cal X})= - \frac{1}{r^{(n-2)/2}} \frac{b_{1}}{\sqrt{Q}}\end{equation}

The normalisation by $1/r^{(n-2)/2}$ means that we do not change $\Psi({\cal
X})$ if we replace ${\cal L}$ by a positive power ${\cal L}^{s}$. Also
this scaling weight is the same as that of the Calabi functional
$$ r^{(n-2)/2} \Vert S(r\omega)\Vert_{L^{2}(r\omega)}
=  \Vert S(\omega)\Vert_{L^{2}(\omega)}. $$
The upshot of this is that in proving our theorem we can always suppose
$r=1$, which we do from now on. By the flatness of the family and the Riemann-Roch
Theorem for $X$, we  can  identify the co-efficient
$a_{0}$ with the volume of $X$ and $a_{1}$ with the integral of $S(\omega)$
for any metric $\omega$ representing $c_{1}(L)$. Thus
$\hat{S}= a_{1}/a_{0}$ and the trivial lower bound (4)  is $\Vert S\Vert\geq
\vert a_{1}\vert/\sqrt{a_{0}}$.

There are two ways in which we can modify a test configuration.
First, we can pull back ${\cal X}$ by a $d$-fold covering of the base. This changes
$A_{k}$ to $dA_{k}$ and plainly does not affect $\Psi$. Second we
can change the $\bC^{*}$-action on the line bundle ${\cal L}$ by a
character $\lambda\mapsto \lambda^{\nu}$ of $\bC^{*}$. This changes
$A_{k}$ to $A_{k} + k \nu 1$ and so
 $b_{i}$ to $\tilde{b_{i}}= b_{i}+ \nu a_{i}$, and
$Q$ to $\tilde{Q}= Q + 2 b_{0}\nu +  a_{0}\nu^{2}$. Thus for the new configuration
$\tilde{{\cal X}}$ we have
\begin{equation}\label{eq:newpsi} \Psi(\tilde{{\cal X}}) = -\frac{b_{1}
 + \nu a_{1}}{\sqrt{Q+ 2 b_{0}\nu  + a_{0}\nu^{2}}}
\end{equation}
We let consider the supremum of this expression over
 $\nu$, initially regarded as a real variable. The analysis of this supremum
 brings in the {\it Futaki invariant} $F_{{\cal X}}$ of the test configuration. This
  is defined to be
 $$     F_{{\cal X}} = b_{1}-  \frac{b_{0} a_{1}}{a_{0}}. $$
 (This terminology differs by a factor of $a_{0}$
  from that in \cite{kn:D2}, where only the sign of $F_{{\cal X}}$ was
  relevant.)
 Then, as the reader will easily verify, if $F_{{\cal X}}< 0$ the expression
 in (9) above is maximised when 
 $$ \nu= \frac{a_{1} Q- b_{0} b_{1}}{a_{0}b_{1}- a_{1} b_{0}}.$$
 By taking a covering we may suppose that this is an integer, so the supremum
 is realised by some test configuration.
 On the other hand if $F_{{\cal X}}\geq 0$ the supremum is not attained for finite
 $\nu$ but occurs in the limit as $\nu$ tends to $\infty$ or $-\infty$, depending on the sign of
 $a_{1}$. In this second case the supremum is just 
 $\vert a_{1}\vert/{\sqrt{a_{0}}}$ which is just the trivial lower bound
 on $\Vert S\Vert$.  In the first case, when the supremum is attained, a
 little calculation shows that the supremum is 
 $$ \left( \frac{\vert a_{1}\vert^{2}}{a_{0}} + 
 \hat{\Psi}_{{\cal X}}^{2}\right)^{1/2}$$
 where we define
 $$  \hat{\Psi}_{{\cal X}}=
 - \frac{F_{{\cal X}}}{\sqrt{Q- b_{1}^{2}/a_{0}}}. $$The denominator here
 has a natural interpretation. We write $\uA_{k}:U_{k}\rightarrow U_{k}$
 for the trace-free part of $A_{k}$, i.e.
 $$ \uA_{k}= A_{k}- \frac{{\rm Tr} A_{k}}{{\rm dim}\ U_{k}} 1. $$
 Then as $k\rightarrow \infty$,
  $${\rm Tr} \uA_{k}^{2}\sim (Q-b_{1}^{2}/a_{0}) k^{n+2}. $$

With this discussion in place we can state our main theorem precisely.
\begin{thm}
If ${\cal X}$ is a test configuration for the pair
$(X,L)$ then for any Kahler metric $\omega$ in the class $c_{1}(L)$
we have
$$  \Vert S(\omega)\Vert_{L^{2}} \geq \Psi({\cal X}). $$

\end{thm}

By the discussion above it is completely equivalent to prove that for any
test configuration
$$  \Vert S(\omega)- \hat{S}\Vert_{L^{2}} \geq -\frac{F_{{\cal X}}}{ Q-\frac{b_{1}^{2}}{a_{0}}},
$$
and it is in this form that we shall prove the result. Notice that
the inequality is vacuous unless $F_{{\cal X}}<0$, that is, unless $(X,L)$
is not K-semistable. We will prove a more general result dealing with
$L^{q}$ norms (loosely analogous to the discussion of other norms in \cite{kn:AB1}).
For an even integer $p$ we define $N_{p}({\cal X})>0$ by the leading term
\begin{equation} \label{eq:pcase} {\rm Tr} \uA_{k}^{p} \sim N_{p}^{p} \ k^{n+p} \end{equation}
in the appropriate Hilbert polynomial. Thus, when $p=2$, 
$$  Q- \frac{b_{1}^{2}}{a_{0}} = N_{2}^{2}. $$
Momentarily re-instating the integer $r$ in the definition of a test configuration,
we define
\begin{equation}\label{eq:hatpsidef} \hat{\Psi}_{p}({\cal X})= -\frac{1}{r^{(n/q)-1}} \frac{
 F_{{\cal X}}}{ N_{p}({\cal X})}. \end{equation}
 Again, the power of $r$ is chosen so that this is unchanged if we replace
 ${\cal L}$ by ${\cal L}^{s}$.
\begin{thm}
If ${\cal X}$ is any test configuration for $(X,L)$ then for any metric
$\omega$ in the class $c_{1}(L)$ we have
$$ \Vert S(\omega)- \hat{S}\Vert_{L^{q}} \geq   \hat{\Psi}_{p}({\cal X}),
$$
where $p$ is any even integer and $q$ is the conjugate exponent to $p$
({\it i.e.} $p^{-1}+q^{-1}=1$ ).
\end{thm}
Just as before, one can check that that the scaling weight of the $L^{q}$
norm of $S(\omega)$ agrees with the power of $r$ in (11), so we can reduce
to the case when $r=1$. By the preceding discussion, Theorem 2 (in the
case when $p=2$) implies Theorem 1, and in the body of the paper below
we prove Theorem 2. (We remark that, when $p\neq 2$, there is no simple
exact relation between the $L^{q}$ norms of $S$ and $S-\hat{S}$. One can
derive slightly different lower bounds for the $L^{q}$ norms of $S$ using
the techniques of this paper, but we find it simpler to work with $S-\hat{S}$.)

If the automomorphism group of the pair $(X,L)$ contains a non-trivial
compact connected subgroup $G$ then there is another relatively elementary
lower bound on the Calabi functional, obtained from the 
 Futaki invariant in its original, differential geometric, form.
 We consider a $G$-invariant metric on $X$. Then we get a Lie algebra homomorphism
 from ${\rm Lie}(G)$ to the functions on $X$, under Poisson bracket.   
Let $\xi$ be an element of ${\rm Lie}(G)$ and let $H$ be the corresponding
Hamiltonian. Then the integrals
$$  \int_{X} S\  H d\mu
\  \ \ \ \  \int_{X} H^{p} d\mu, $$
are topological invariants of the data: they do not change as we vary
$\omega$ among $G$-invariant metrics in the same Kahler class.
   By Holder's inequality we have
   \begin{equation}\label{eq:Holderbound} \Vert S-\hat{S}\Vert_{L^{q}} \geq \frac{1}{\Vert H-\hat{H}\Vert_{L^{p}}}
   \int_{X} (S-\hat{S}) H d\mu \end{equation}
   where $\hat{H}$ is the average value of $H$, and the right hand side
   is a topological invariant of the data.
    
We will now explain how to derive this as a special case of our
 Theorem 2. By a density argument we can suppose that $\xi$ generates a circle
 action on $(X,L)$ which we extend to a holomorphic 
 $\bC^{*}$ action. Then we define a test configuration by taking
 the product ${\cal X}= X\times \bC$ but using this non-trivial
 $\bC^{*}$ action. Thus in this case the central fibre $(X_{0}, {\cal L})$ is just the
 original pair $(X,L)$, with the given $\bC^{*}$ action. The key point
 now is that 
 $$ b_{0}= \int_{X} H d\mu\ , \ b_{1}= \int_{X} S H d\mu, $$
 so that
 $$   F_{{\cal X}} = \int_{X} (S-\hat{S}) H d\mu. $$
  Thus the algebro-geometric definition of the Futaki
 invariant reduces to the differential geometric one.  This is explained
  in Prop. 2.2.2 of \cite{kn:D2} and we will obtain a generalisation
 of the fact in Proposition 3 below. Moreover, the same discussion shows that, for
 even integers $p$,
 $$  N_{p}= \Vert H- \hat{ H} \Vert_{L^{p}}. $$
Changing the sign of $\xi$, if necessary, we can suppose that the
Futaki invariant is negative. Then 
 (12) follows  as a special case of our Theorem 2. In fact our Theorem gives rather more even
in this case since the lower bound is obtained for arbitrary metrics in
the Kahler class, not just the $G$-invariant ones.

We return now to situation considered by Atiyah and Bott and explain how
their result can be cast in a similar form to ours. Let
${\cal F}$ be a flag of sub-bundles in $E\rightarrow \Sigma$, as before.
Let $W=(w_{1}, \dots, w_{q})$ be  a vector of strictly {\it increasing}
 integers. Then one can construct a degeneration ${\cal E}= {\cal E}({\cal
 F}, W)$ from this data. This
  is a $\bC^{*}$-equivariant bundle over $\bC^{*}\times \Sigma$ 
 which is isomorphic  to $E$ on each slice $\{ t \} \times \Sigma$ for
 $t\neq 0$ and to 
 $E_{0}=Q_{1} \oplus Q_{2} \dots \oplus Q_{q}$ when $t=0$, with the
 property that the 
 $\bC^{*}$-action on $E_{0}$ has weight $w_{i}$ on the summand
 $Q_{i}$. We fix a square root $K^{1/2}_{\Sigma}$ of the canonical bundle
 and any line bundle $L\rightarrow \Sigma$ of degree $1$.
 Now we consider the $\bC^{*}$ action on the vector spaces $$ U_{k}= H^{0} (E_{0} \otimes K^{1/2}_{\Sigma}\otimes
 L^{k}), $$ with generators $A_{k}:U_{k}\rightarrow U_{k}$. Following just
 the same pattern as before, we look at the large $k$ behaviour,
 
$$   {\rm dim}\ U_{k}= \alpha_{0} k+ \alpha_{1}, $$ 
$$  {\rm Tr} A_{k}= \beta_{0} k + \beta_{1}, $$
$$ {\rm Tr} A_{k}^{2}\sim Q k, $$
and we define $${\Psi}({\cal F}, W)= -\frac{\beta_{1}}{\sqrt{Q}}. $$
\begin{lem}
For any bundle $E\rightarrow \Sigma$ we have
$$  \sup_{{\cal F}, W} \Psi({\cal F}, W)= \sup_{{\cal F}, {\rm dec}} \Phi({\cal
F}), $$
where on the left hand side the supremum is taken over arbitrary flags
${\cal F}$
and increasing weight vectors $W$ and on the right hand side the supremum
is taken over slope-decreasing flags ${\cal F}$.  
\end{lem}
Given a flag ${\cal F}$ and weight vector $W$ we have
$$  U_{k}= \bigoplus H^{0}( Q_{i} \otimes K^{1/2}_{\Sigma} \otimes L^{k}),
$$
with the action of weight $w_{i}$ on the the ith summand.  
It follows that 
$$ \alpha_{0}= \sum r_{i}, \alpha_{1} = \sum d_{i}, \beta_{0}= \sum r_{i}
w_{i}, \beta_{1}= \sum d_{i} w_{i}, Q= \sum r_{i} w_{i}^{2}, $$
 where $Q_{i}$ has rank $r_{i}$ and degree $d_{i}$.
  So we obtain
\begin{equation}\label{eq:psiform} \Psi({\cal F}, W) = - \frac{\sum d_{i} w_{i}}
{\sqrt{\sum r_{i} w_{i}^{2}}}. \end{equation}
Now suppose that ${\cal F}$ is slope-decreasing. We maximise
$\Psi({\cal F}, W)$ over weights $w_{i}$, initially regarded
as arbitrary real numbers.  The maximum occurs at $w_{i}^{{\rm max}}= -C d_{i}/r_{i}$
for any constant $C$. Taking a suitable $C$ we can suppose  that
$w_{i}^{{\rm max}}$ are integers, but more crucially the condition that
${\cal F}$ was a slope-decreasing flag means that the $w_{i}^{{\rm max}}$
are increasing. Then we have
$$ \Psi({\cal F}, W^{{\rm max}}) = \Phi({\cal F}),$$
so we have established one half of the Proposition, {\it i.e} 
$\sup \Psi \geq \sup \Phi $. We prove the other half by induction
on the length $q$ of the flag. We take as inductive proposition the statement
that for any flag ${\cal F}$ of length $q$ and any {\it weakly} increasing weight
vector $W$ ({\it i.e} with $w_{1}\leq \dots \leq w_{q})$ we have
$\Psi({\cal F},W) \leq \sup_{{\cal F}, {\rm dec}} \Phi$, where
$\Psi({\cal F}, W)$ is defined by (13). This is clearly true for $q=1$.
 Suppose that ${\cal F}$ is any flag and
$W^{0}$ is a weakly-increasing weight vector. We  maximise the expression for $\Psi({\cal F}, W)$
over all weakly-increasing weight vectors. If ${\cal F}$ is slope-decreasing
then we are in the same position as above and the maximum realises
$\Phi({\cal F})$; so in this case
$$ \Psi({\cal F},W^{0}) \leq \Psi({\cal F}, W^{max})\leq \sup \Phi$$
as desired.  If ${\cal F}$ is not slope increasing then the
maximum occurs at a vector $W^{max}$ with $w^{max}_{i}=w^{max}_{i+1}$ for
some $i$. Then let ${\cal F}'$ be the flag obtained from ${\cal F}$ by deleting
$E_{i}$ and let $W'$ be the weight vector for ${\cal F}'$ obtained from
$W^{max}$ in the obvious way. In the associated sum of quotients we are
replacing the original pair $Q_{i}\oplus Q_{i+1}$ by a bundle
$Q'_{i}= E_{i+1}/E_{i-1}$. The rank of $Q'_{i}$ is the sum $r_{i}+r_{i-1}$
and the degree of $Q'_{i}$ is $d_{i}+d_{i-1}$. Thus
$\Psi({\cal F}', W') = \Psi({\cal F}, W^{max})$. So we have
$$\Psi({\cal F}, W^{0}) \leq \Psi({\cal F}, W^{max})= \Psi({\cal F}', W')
\leq \sup \Phi, $$
where in the last inequality we use the inductive hypothesis. This completes
the proof.

 In sum, the version of the 
Atiyah-Bott result which is closer to our Theorem is the statement that
$$  \inf_{A} \Vert F(A)\Vert = \sup_{{\cal F}, W} \Psi({\cal
F}, W), $$
which is  completely equivalent to the previous formulation by the elementary
Lemma above.  
One could probably prove one direction here (i.e. that
$\inf_{A} \geq \sup_{{\cal F}, W}$) by using a finite-dimensional
approximation
argument in the manner of Wang \cite{kn:W}, just as we will do for the
Calabi functional. However there would not be much point to this,
in view of the simple and direct differential geometric proof.

\section{Differential Geometric  asymptotics}

In this section we relate the Calabi functional (and its $L^{q}$ variants)
 to the norm of a matrix
$\uM$ associated to a projective variety. Let $z_{\alpha}$ be standard
homogeneous coordinates on $\bC\bP^{N}$ and let
let $h_{\alpha \beta}$ be the function
$$  h_{\alpha \beta} = \frac{z_{\alpha} \overline{z}_{\beta}}{\Vert \underline{z}\Vert^{2}}$$
on $\bC\bP^{N}$. For a smooth projective variety $V\subset \bC\bP^{N}$
of dimension $n$
 we define a self-adjoint matrix $M=M(V)$ with entries
$$   M_{\alpha \beta} = \int_{V} h_{\alpha \beta} d\mu_{FS}, $$
where $d\mu_{FS}$ is the standard volume form induced from the Fubini-Study metric
(normalised so that the volume of $V$ is equal to its degree divided by
$n!$). Let $\uM$ be the trace-free part of $M$. Since the sum of the
$h_{\alpha \alpha}$ is $1$ we have
$\uM = M - \frac{{\rm Vol}(V)}{N+1} 1$.
We recall that for any $q>1$ the $q$-norm of a  self-adjoint matrix $T$
is defined by 
$$  \Vert T \Vert_{q}^{q}= \sum \vert \lambda_{\alpha}\vert ^{q} , $$
where $\lambda_{\alpha}$
are the eigenvalues, repeated according to multiplicity.

Now let $(X,L)$ be an abstract polarised variety, as before. For large $k$
the sections of $L^{k}$ define a projective embedding of $X$,  i.e a choice
of basis of $H^{0}(L^{k})$ yields a specific projective variety
$V$. 
\begin{prop}
Let $\omega_{0}$ be a Kahler metric on $X$ in the class $c_{1}(L)$.
Then for large enough $k$ there is a basis of $H^{0}(L^{k})$
 yielding a projective embedding $X\rightarrow V_{k}\subset \bC\bP^{N_{k}}$
 with for, any $q\geq 1$,
 $$ \Vert \uM(V_{k})\Vert_{q} \leq k^{(n/q)-1} \Vert
  S(\omega_{0})- \hat{S}\Vert_{L^{q}}
  + O(k^{(n/q)-2}).
 $$
\end{prop}

The proof is a straightforward application of the
  asymptotic expansion for the density of states function (see references
  in Section 1)
which we now recall. Let $\vert \ \vert_{0}$ be a Hermitian metric on
$L$ whose associated curvature form is $-2\pi i \omega_{0}$. We write
$\vert \ \vert_{0}$ also for the induced metric on $L^{k}$. Then endow 
 $H^{0}(L^{k})$ with the standard $L^{2}$-norm
defined by the volume form $d\mu_{0}$ of the fixed metric $\omega_{0}$
and the fibre metric $\vert \ \vert_{0}$ on $L^{k}$.
 Define a function
$\rho_{k}$ on $X$ by
$$  \rho_{k}= \sum \vert s_{\alpha}\vert^{2}, $$
for any orthonormal basis $s_{\alpha}$ of $H^{0}(L^{k})$.
Then the statement we need is that
$$\rho_{k} = k^{n} ( 1+ k^{-1} \eta_{k})$$
where the functions $\eta_{k}$ converge (in $C^{\infty}$) to 
$ S(\omega_{0}) $ as $k\rightarrow \infty$; in fact
$\eta_{k}=S+ O(k^{-1})$. 

The construction of the projective embedding we need is the most
 obvious one.  For each $k$, we take any orthonormal basis $s_{\alpha}$ of 
 $H^{0}(X,L^{k})$ and define a new fibre metric on the line bundle $L$ by
 $$  \vert \ \vert_{*}^{2}= \frac{1}{\rho_{k}^{1/k}}. $$
 Then, denoting the induced metric on $L^{k}$ also by $\vert \ \vert_{*}$,
  $$ \sum \vert s_{\alpha}\vert_{*}^{2} = 1. $$
  The pull-back of the Fubini-Study metric under the embedding is $k\omega$ where
 $$   \omega= \omega_{0} + k^{-1} \frac{i}{2\pi}
  \overline{\partial}\partial (\log (1+ k^{-1}\eta_{k})). $$
 Write the volume form of $\omega$ as
 $$ (1+ k^{-2} \nu_{k}) d\mu_{0}. $$
 It is clear, from the fact that the sequence $\eta_{k}$ is bounded, that
 the $\nu_{k}$ are bounded. The functions $h_{\alpha \beta}$ pull back to
 under the projective embedding to $(s_{\alpha},s_{\beta})_{*}$. So we
 have
 $$\int_{V_{k}}h_{\alpha \beta} d\mu_{FS} =  \int_{X} (s_{\alpha},s_{\beta})_{0}
 \left(\frac{1+ k^{-2} \nu_{k}}{1+ k^{-1} \eta_{k}}\right) d\mu_{0}.
   $$
   We do not change the norm of $\uM(V_{k})$ if we apply a unitary transformation
   to $\bC^{N+1}$, in other words if we make a different choice of 
   orthonormal
   basis $s_{\alpha}$. Thus we can choose the  basis so that $M$ is a diagonal
   matrix with diagonal entries $m_{\alpha}$ where
   $$ m_{\alpha} =  \int_{X} \vert s_{\alpha} \vert_{0}^{2}
    \left( \frac{1+
   k^{-2} \nu_{k}}{1+ k^{-1} \eta_{k}}\right) d\mu_{0}.$$
   
   Now the dimension $N+1$ is the integral of $\rho_{k}$ over $X$, and it
   follows that the trace-free part $\uM$ of $M$ is the diagonal matrix
   with diagonal entries
   $$  \um_{\alpha} = m_{\alpha} - \frac{1}{1+ \hat{\eta_{k}} k^{-1}}, $$
   where $\hat{\eta}_{k}$ is the average value of $\eta_{k}$.
   Then we obtain
   $  \um_{\alpha} = k^{-1} (b_{\alpha} + \epsilon_{\alpha})$ where
   $$ b_{\alpha} = \int_{X} \vert s_{\alpha}\vert_{0}^{2} (\hat{\eta_{k}}-
   \eta_{k})
   d\mu_{0}, $$
   and $\epsilon_{\alpha}$ is $O(k^{-1})$. 
   Now write
   $$  \vert s_{\alpha}\vert_{0}^{2}\  \vert \eta_{k} - \hat{\eta_{k}}\vert=
     \vert s_{\alpha}\vert_{0}^{2/p} \vert s_{\alpha}\vert_{0}^{2/q} \vert
     \eta_{k}- \hat{\eta_{k}}\vert, $$
   where $p$ is the index conjugate to $q$. Apply
      Holder's inequality to get
     $$   \vert b_{\alpha}\vert \leq
\left( \int_{X} \vert s_{\alpha}\vert_{0}^{2} d\mu_{0}\right)^{1/p} \left(
\int_{X} \vert s_{\alpha}\vert_{0}^{2} \vert \eta_{k} - \hat{\eta_{k}}\vert^{q}
d\mu_{0}\right)^{1/q}.$$
But since the $s_{\alpha}$ are orthonormal this gives
$$ \vert b_{\alpha}\vert^{q} \leq \int_{X} \vert s_{\alpha}\vert_{0}^{2}
\vert \eta_{k}-\hat{\eta_{k}} \vert^{q} d\mu_{0}.$$
Summing over $\alpha$ and using the asymptotic statement again in the weak
form
$$ \sum \vert s_{\alpha}\vert_{0}^{2} = k^{n} + O(k^{n-1}), $$ we obtain
$$ \sum_{\alpha} \vert b_{\alpha}\vert^{q} \leq (k^{n}+ C k^{n-1}) \int_{X} \vert
\eta_{k} - \hat{\eta_{k}}\vert^{q}   d\mu_{0} . $$
Now $$\Vert \uM\Vert_{q} \leq k^{-1} ( \Vert B\Vert_{q} + \Vert E \Vert_{q})$$
where $B,E$ are the diagonal matrices with entries $b_{\alpha}, \epsilon_{\alpha}$
respectively. We have
$$ \Vert B \Vert_{q} \leq k^{n/q} (1+ C k^{-1})^{1/q} \Vert \eta_{k}-\hat{\eta_{k}}\Vert_{L^{q}}
\leq k^{n/q} (1+ \frac{C}{q} k^{-1}) \Vert \eta_{k}-\hat{\eta_{k}}\Vert_{L^{q}},
$$
and $\Vert E \Vert = O(k^{(n/q)-1})$ since the dimension $N+1$ is
$O(k^{n})$. This gives
$$ \Vert \uM \Vert_{q} \leq k^{(n/q)-1} \Vert \eta_{k}-\hat{\eta_{k}}\Vert_{L^{q}}
+ O(k^{(n/q)-2}), $$
and our result follows from the fact that $\eta_{k}-\hat{\eta_{k}}= (S(\omega_{0})-\hat{S})+
O(k^{-1}) $.

    \section{The finite-dimensional argument}

    Consider a $\bC^{*}$ action on $\bC\bP^{N}$ induced by a $1$-parameter
    subgroup $\rho:\bC^{*}\rightarrow GL(N+1)$. Let $V$ be a smooth
    $n$ dimensional projective
    variety and set $V^{t} = \rho(t)(V)$. Then it follows from standard
    theory that the $ V^{t}$ converge as $t\rightarrow 0$ to some 
    algebraic cycle $V^{0}$. Thus $V^{0}$ is a formal sum
    $$ V^{0} = \sum m_{i} W_{i} $$
    where the mulptiplicities $m_{i}$ are positive integers and
    $W_{i}$ are irreducible $n$-dimensional projective varieties. This convergence can
    be understood at various levels, but the crucial point for us is that
    the $V^{t}$ converge to $V^{0}$ in the sense of currents. So for any
    smooth test form $\phi$
    $$ \int_{V^{t}} \phi \rightarrow \int_{V^{0}}\phi = \sum m_{i} \int_{W_{i}
   } \phi. $$
    Now suppose that $\rho$ maps $S^{1}\subset \bC^{*}$ to the unitary
    group. Thus the infinitesimal generator $A$ of $\rho$ is a Hermitian
    matrix; as usual, we write $\uA$ for the trace-free part of $A$. Let $h$ be the real valued function
    $$  h=\sum_{\alpha \beta} A_{\alpha \beta} h_{\alpha \beta} $$
    on $\bC\bP^{N}$. This is a Hamiltonian for the action of $S^{1}$ on
    the projective space. Let 
    \begin{equation}\label{eq:Idef}  I(A,V^{0}) = \int_{V(0)} h\ d\mu_{FS} =\frac{1}{n!}
     \int_{V(0)} h \ \omega_{FS}^{n}, \end{equation}
     where $\omega_{FS}$ is the Fubini-Study form, and set
     $$FCh(A,V^{0}) =  \frac{{\rm Vol}(V^{0})}{N+1} {\rm Tr}(A) - I(A,V^{0}), $$
      where the volume of $V^{0}$ has the obvious meaning.  
     \begin{prop}
     Suppose that $FCh(A,V^{0})<0$. Then for any conjugate indices $p,q$
     we have
     $$\Vert \uM(V)\Vert_{q} \geq \frac{- FCh(A, V(0))}{\Vert \uA \Vert_{p}} $$     
     \end{prop}

     To prove this we consider the function
     $$   f(t)= {\rm Tr}( \uA\  \uM(V^{t} )), $$
     for $t\in {\bf R}^{*}$. The crucial point is that this is {\it 
     increasing} with $t$.
     From one point of view this follows from the fact that $\uM$ is
     a moment map, and from the general theory of such maps, see Section
     6.5.2 in \cite{kn:DK},
     for example. If $\mu$ is a moment map for an isometric action of a
     group $G$ on a Kahler manifold we have, for any $a\in {\rm Lie}(G)$,
     $$   \frac{d}{ds} \langle \mu( {\rm exp}(sa)X), a \rangle= \vert \frac{d}{ds}
     {\rm exp}(sa) X \vert^{2} \geq 0.$$
     Then our assertion follows (taking $t=e^{s}$) from the fact established
     in \cite{kn:D3} that
     $\uM$ is a moment map for the $SU(N+1)$ action on the set of varieties
     projectively
     equivalent to $V$, with a suitable Kahler structure. 
     An equivalent fact was used in the earlier work of Zhang \cite{kn:Z},
     and the relation between this and the constructions of \cite{kn:D3}
     is explained in \cite{kn:PS1}, \cite{kn:PS2}.  This monotonicity property
     is often stated in the literature 
     as the {\it convexity} of a certain function, and other direct proofs
     are given in \cite{kn:D4}, Prop.1,  and \cite{kn:PS2}, Lemma 3.1, 
     so we will not discuss the matter further here.
     Given this monotonicity property we argue as follows. We have
     $${\rm Tr} (A M(V^{t}))= \int_{V(e^{t})} h_{A} d\mu_{FS} $$
     and $\uM= M- \frac{{\rm Vol}(V)}{N+1} 1$. So
     $$ f(t)= {\rm Tr} (\uA\ \uM(V^{t}))= {\rm Tr}(A \uM(V^{t}))=
     \int_{V^{t}} h_{A} d\mu_{FS} - \frac{{\rm Vol}}{N+1} {\rm Tr} A.
     $$
     Hence the limit of $f(t)$ as $t\rightarrow 0$ is
          $-FCh(A, V^{0})$, which is positive by hypothesis.
      Since the function $f$ is increasing 
      $$ \vert f(t) \vert = f(t) \geq - FCh(A,V^{0})$$
      for all $t$. In particular, taking $t=1$, we have
      $$  {\rm Tr}(\uA \ \uM(V)) \geq -FCh(A, V^{0}). $$
      Now use the fact that for, any Hermitian matrices $S,T$
      $$  \vert {\rm Tr} (ST)\vert \leq \Vert S\Vert_{p} \Vert T \Vert_{q},
      $$ to obtain the required result.

     We remark that all of this dicussion can be placed in the context
     of the action of $GL(N+1)$ on the space of Chow vectors, as explained
     in \cite{kn:PS1}, \cite{kn:PS2}, and the criterion $FCh(A,V^{0})<
     0$ is the standard Hilbert-Mumford criterion for a destabilising $1$-parameter
     subgroup. On the other hand the criterion can also be placed in
      a dynamical,
     Riemannian geometry context. Let $P$ be a compact Riemannian manifold,
     $S\subset P$ a submanifold and let $h$ be a real-valued function on
     $P$. Let $\phi_{s}$ be the flow generated by the ($h$-decreasing)
     gradient of $h$ and suppose that the $\phi_{s}(S)$ converge 
     to some varifold $S'$ as $s\rightarrow \infty$. Say that $S$ has property
     (*) with respect to $h$ if the mean value of $h$ over $S'$ is less
     than the mean value of $h$ on $P$. Then a variety $V\subset \bC\bP^{N}$
     is Chow stable if it has property (*) with respect to all the functions
     $h_{A}$. This is because
     $$   -FCh(A,V^{0})= {\rm Vol}(V^{0}) \left( \frac{I}{{\rm Vol}(V^{0})}
     - \frac{{\rm Tr}(A)}{N+1}\right), $$
     and ${\rm Tr}(A)/ N+1$ is the mean value of $h_{A}$ on $\bC\bP^{N}$.

     \section{Algebro-geometric asymptotics}
     
     We now go back to a test configuration ${\cal X}$ for $(X,L)$.
      We can
     suppose that the parameter $r$ is $1$. The essential point 
      is that for large enough $k$ the configuration can be realised by
      a $\bC^{*}$-action on an ambient projective space. That is to say,
      ${\cal X}$ can be embedded as  a $\bC^{*}$-invariant subscheme in the product
      $\bP(U_{k}^{*}) \times \bC$, extending the embedding
      of  the central fibre $X_{0}\subset
      {\cal X}$  by the complete linear system $U_{k}$ in 
      $\bP(U_{k}^{*})=\bP(U_{k}^{*}) \times \{0\}. $
      This is  explained in {\cite{kn:RT}. In essence we consider the
      $\bC^{*}$-equivariant bundle ${\cal U}=\pi_{*} {\cal L}^{k}$ over $\bC$
      and pick an equivariant trivialisation ${\cal U}\cong \bC \times U_{k}$.
       There is no loss
      of generality in supposing that ${\cal L}^{k}$ is very ample on all fibres
      for all $k\geq 1$; then we obtain the embedding of ${\cal X}$ from
      the fibrewise embeddings of $X_{t}$ in $\bP({\cal U}_{t}^{*})$ under this
      trivialisation. Let
      $\vert X_{0}\vert$ denote the cycle in $\bP(U_{k}^{*})$ associated to the scheme
      $X_{0}$. Then we are precisely in the situation considered above
      with a $1$-parameter subgroup acting on $\bP^{N+1}$, generated by $A_{k}$.
       If we write $V^{t}$ for
      the image of $X_{t}\times \{t\}\subset {\cal X} \subset \bP \times
      \bC$ under the projection map to $\bP$ then $V^{t}= \rho(t)(V_{1})$
      and the $V^{t}$ converge to the cycle $\vert X_{0}\vert$ as $t\rightarrow
      0$.

      We need a general fact about the choice of equivariant trivialisation
      of the bundle ${\cal U}$. 
      \begin{lem}
      Let $E$ be $\bC^{*}$-equivariant bundle over $\bC$ and let $H$ be
      a Hermitian metric on the fibre $E_{1}$. Then there is an equivariant
      trivialisation $E\cong \bC \times E_{0}$ which takes $H$ to a Hermitian
      metric on the central fibre $E_{0}$ which is preserved by the action of $S^{1}\subset
      \bC^{*}$ on $E_{0}$.
      \end{lem}
       To see this consider the weight spaces
       $E_{0}= \bigoplus V_{i}$ say, where $\bC^{*}$ acts with weight
       $w_{i}$ on $V_{i}$ and the ordering is chosen so that
       $w_{1}< w_{2} < \dots$. Let ${\cal F}$ be the flag
       $$  V_{1} \subset V_{1} \oplus V_{2} \subset V_{1}\oplus V_{2} \oplus
       V_{3} \dots, $$
       in $E_{0}$ and let $\alpha:E_{0}\rightarrow E_{0}$
       be a linear map which preserves the flag ${\cal F}$. Thus
       $\alpha$ has a block matrix description $\alpha_{ij}: V_{i} \rightarrow
       V_{j}$ with $\alpha_{ij}=0$ for $i>j$. For any $t\in \bC$ we define
       $\alpha(t)$ with blocks $t^{w_{j}-w_{i}} \alpha_{ij}$.
       This defines a $\bC^{*}$ equivariant automorphism of the trivial
       bundle $\bC \times E_{0}$ equal to $\alpha$ on the fibre $\{1\}\times
       E_{0}$. Conversely, all equivariant automorphisms arise in this
       way. What this means is that in the fibre $E_{1}$ 
       of our equivariant bundle $E$ there is a canonical flag and a choice
       of equivariant trivialisation is equivalent to a choice of compatible
       direct sum splitting of $E_{1}$. Thus the proof of the Lemma is
       simply to take the direct sum splitting furnished by the succesive
       orthogonal complements in the flag in $E_{1}$ using the metric $H$,
       and then take the corresponding equivariant trivialisation.

Applying the Lemma to the bundle ${\cal U}$ we see that, given any metric
on $H^{0}(X,L^{k})$, we may choose our representation of ${\cal X}$ so that
the $1$-parameter subgroup takes the circle to unitary transformations
of $U_{k}$, with respect to the metric arising from 
 the identification of $U_{k}$ with $H^{0}(X,L^{k})$. (Of course we can
  choose a basis so that the given metric is identified with the standard
 one on $\bC^{N+1}$.)      
       Then  we have a numerical invariant
       $$ I_{k}= I(A_{k},\vert X_{0}\vert)=\int_{\vert X_{0}\vert} h d\mu_{FS}, $$
       as above.
       \begin{prop}
       The integral $I_{k}$ is equal to the leading term $b_{1} k^{n+1}$
        in the Hilbert polynomial for ${\rm Tr} (A_{k})$.
       \end{prop}
       Assuming this for the moment we go on to complete the proof of our
       main result. Observe that it is equivalent to consider a fixed embedding
       $X\rightarrow V=V^{1}$ with a variable metric on the underlying
       vector space $\bC^{N+1}$ (as we are doing here) or to consider a
       fixed metric on $\bC^{N+1}$ and varying the embedding by projective
       transformations (as we considered in Section 3). 
        Set $$FCh_{k}= FCh(A_{k}, \vert X_{0}\vert)=  \frac{{\rm Vol}(X,k\omega)}{{\rm
       dim} U_{k}}{\rm Tr} A_{k}-I_{k},$$ as in the previous section. 
       Thus
       $$ FCh_{k}=  \left( \frac{a_{0} k^{n}}{a_{0} k^{n}
       + a_{1} k^{n-1} + \dots} \right) ( b_{0} k^{n+1} + b_{1} k^{n} +
       \dots)- b_{0}  k^{n+1}.  $$
       So
       \begin{equation}\label{eq:Fconv}  FCh_{k}= (b_{1} - \frac{a_{1} b_{0}}{a_{0}}) k^{n} + O(k^{n-1})
       = F_{{\cal X}} k^{n} + O(k^{n-1}). \end{equation}
       Now Theorem 2 is vacuous if $F_{{\cal X}}\geq 0$, so suppose that
       $F_{{\cal X}}<0$. Then (15) means that $FCh_{k}<0$ for large $k$. Thus
       we can apply Proposition 2 to deduce that for the embedding $X\rightarrow V=V(1)$ 
       $$ \Vert \uM(V) \Vert_{q} \geq - \frac{FCh_{k}}{ \Vert \uA_{k}\Vert_{p}}.
       $$
       If $p$ is an even integer then $\Vert A_{k}\Vert_{p}^{p} = {\rm Tr}
       \uA^{p}$, so 
       $$  \frac{FCh_{k}}{ \Vert \uA_{k}\Vert_{p}} = k^{(n/q)-1}
        \frac{F_{{\cal X}}}{N_{p}}
         + O(k^{n/q)-2}). $$
         
        Thus \begin{equation} \label{eq:qbound}\Vert \uM(V)\Vert_{q} \geq - k^{(n/q)-1} \frac{F_{{\cal
        X}}}{N_{p}} + O(k^{(n/q)-2}).\end{equation}
        Then Theorem 2 follows from this lower bound combined with Proposition
        1, since if there was a Kahler metric $\omega_{0}$ with 
        $ \Vert S-\hat{S}\Vert_{L^{q}} < \hat{\Psi}_{p}({\cal X})$ we would
        get, for large $k$,  an embedding $V$ with $\Vert \uM(V)\Vert_{q}< - k^{(n/q)-1}
        F_{{\cal X}}/ N_{p}$ in contradiction to (16).

       \subsection{Proof of Proposition 3}

        First observe that if we have proved that ${\rm Tr}(A_{k}) \sim
        k^{n+1} I_{1}$ then replacing $L$ by $L^{s}$ it will follow that
        ${\rm Tr}(A_{sk})\sim (sk)^{n+1} I_{s}$ and so that $I_{s}= s^{n+1}
        I_{1}$. (This can also be seen directly. If we asume that the powers
        of sections in $H^{0}(X_{0}, {\cal L})$ generate $H^{0}(X_{0}, {\cal
        L}^{s})$ then a choice of metric on the first space yields a natural
        metric on the second. With these metrics the integrand defining $I_{s}$ is 
        $s^{n+1}$ times
        that defining $I_{1}$, pointwise on $\vert X\vert$.)
        To sum up, our task is to establish the asymptotic relation
        \begin{equation} \label{eq:task}  {\rm Tr}(A_{k}) \sim k^{n+1} I, \end{equation}
        where we recall that \begin{itemize}
        \item $X_{0}$ is a $\bC^{*}$-equivariant subscheme of $\bC\bP^{N}$
        for the action generated by $A=A_{1}$ on $\bC^{N+1}= H^{0}(X_{0},
        {\cal O}(1))$;
        \item $A_{k}$ is the generator of the induced action on $H^{0}(X_{0},
        {\cal O}(k)) $;
        \item We fix an $S^{1}$-invariant metric on $\bC^{N+1}$. Then
         $I$ is the integral of $h\omega_{FS}^{n}/n!$ over the cycle
        $\vert X_{0}\vert$ associated to $X_{0}$, where $\omega_{FS}$ is
        the Fubini-Study metric on $\bC\bP^{N}$ and $h$ is the function
        on $\bC\bP^{N}$
        associated to $A$.
        \end{itemize}
        
        If $X_{0}$ is smooth the relation (17) is rather standard. The key
        point is that the function $h$ is a Hamiltonian for the $S^{1}$
        action on $\bC\bP^{N}$. The desired result is then obtained from
        the equivariant Riemann-Roch Theorem and the de Rham model for
        equivariant cohmology, as explained by Atiyah and Bott in 
        \cite{kn:AB2};
        see the discussion in Section 2 of \cite{kn:D2}. 
        Our problem is to extend this discussion to the case when $X_{0}$
        is a singular variety or scheme. What we do take as known is the
        corresponding non-equivariant result. That is, if $Z\subset \bC\bP^{M}$
        is an $m$-dimensional projective scheme then 
        \begin{equation} \label{eq:known}   {\rm dim}\ H^{0}(Z, {\cal O}(k)) \sim \frac{D}{m!} k^{m}, 
        \end{equation}
        where $D$ is the degree of $Z$, which is the integral over the
        cycle $\vert Z\vert$ of $\omega_{FS}^{m}$. (The proof of this can
        be reduced to the case $m=0$ by taking hyperplane sections. Thus
        the assertion is essentially that the notions of {\it multiplicity}
         defined
        algebraically or by currents agree.)
        We will explain how to
         prove the equivariant result by reducing to (18).
        
        We consider the following general situation. Let $P\rightarrow
        B, Q\rightarrow F$ be a pair of principle $S^{1}$-bundles over 
         manifolds $B, F$. So we have $S^{1}$ actions $\sigma_{P}, \sigma_{Q}$
         say and vector fields $v_{P}, v_{Q}$ on $P,Q$. Suppose in addition
         that
         $Q\rightarrow F$ is an $S^{1}$-equivariant circle bundle, so we
         have another action $\rho$ on $Q$, commuting with $\sigma_{Q}$,
         and another vector field, $w$ say, on $Q$. We also denote the induced
         action on $F$ by $\rho$. Now we can form the
         associated bundle
         $$   M=P\times_{\sigma_{P},\rho} F. $$
         Thus $M$ is a bundle over $B$ with fibre $F$. We can also form
         $$ \Pi= P\times_{\sigma_{P}, \rho} Q. $$
         Then $\Pi$ is a circle bundle over $M$.  Suppose now that
         we have connections on the bundles $P\rightarrow B, Q\rightarrow
         F$ and that the connection on $Q$ is preserved by $\rho$. Thus
         we have $1$-forms $\alpha_{P}, \alpha_{Q}$ on $P,Q$ and
         $$  \alpha_{P}(v_{P})=1, \alpha_{Q}(v_{Q})=1, \alpha_{Q}(w)= -H,
         $$
         where $H$ is an $S^{1}\times S^{1}$-invariant function on $Q$,
         which can also be regarded as a function on $F$ or on $M$. Now
         the connection on $P$ defines a splitting of the
         tangent bundle of  the associated bundle  $M$, which we express
         rather loosely as
         $$ TM= TF \oplus TB. $$
         Via this decomposition, the curvature forms $\omega_{P}, \omega_{Q}$ can naturally be
         regarded as $2$-forms on $M$. The key point is that there is a
         connection on $\Pi\rightarrow M$ with curvature
         $$ \Omega= \omega_{Q} + H \omega_{P}. $$
         To see this we pull everything back to $P\times Q$ where we
          consider
         the $1$-form
         $$  \beta= \alpha_{Q} + H \alpha_{P}, $$
         (making an obvious simplification in notation). Then $\beta$ vanishes
         on the generator $v_{P}+ w$ of the action $(\sigma_{P},\rho)$
         and so  descends to a form $\underline{\beta}$ on
         $\Pi=P\times_{\sigma_{P},\rho} Q$.
         Since $\beta(v_{Q})=1$ the equivariant $1$-form $\underline{\beta}$ furnishes
         a  connection on $\Pi$. We have (again making various abuses of
         notation)
         $$  d\beta= \omega_{Q} + H \omega_{P} + dH \wedge \alpha_{P}.
         $$                
         This gives the lift of the curvature form of the connection on
         $\Pi$ to $P\times Q$. The definition of the horizontal subspace
         defining the splitting $TM=TF\oplus TB$ means that $\alpha_{P}$
         vanishes on the vectors representing $TM$, so we see that the connection
         $\underline{\beta}$ has curvature $\Omega$. 
         
         Now suppose that  $B$ is a compact oriented $2r$-manifold and let $S$ be an $S^{1}$-invariant
   oriented         $2n$-dimensional submanifold of $F$, not necessarily closed. This defines a corresponding
          submanifold $\tilde{S}$ in $M$ (so $\tilde{S}$ is a bundle over
          $B$ with fibre $S$). Then it is clear, by integrating over the
          fibres, that
          \begin{equation}\label{eq:intid} \frac{1}{(n+r)!}
          \int_{\tilde{S}} \Omega^{n+r} = \frac{1}{r!}\int_{B} 
          \omega_{P}^{r} \
          \ 
          \frac{1}{n!}\int_{S} H^{n} \omega_{Q}^{n}. \end{equation}
          
          We apply this to the case when $B$ is the Riemann sphere
          $\bC\bP^{1}$ and $F$ is $\bC\bP^{N}$ with the bundles $P, Q$
          being the Hopf fibrations.
          We take the action $\rho$ to be that induced by the given action
          on $\bC^{n+1}$. Then the function $H$ considered above becomes
          the function $h$ by the standard discussion of Hamiltonians and
          equivariant classes, as in \cite{kn:AB2}. We then have a manifold
          $M$, which is just the projectivization of a vector bundle over
          $\bC\bP^{1}$, and a circle bundle $\Pi\rightarrow M$ which clearly
          corresponds to a holomorphic line bundle say ${\cal V}\rightarrow
          M$. Notice that, going back to our original data, we could change
          the generator $A$ to $A+\nu 1$ for any integer $\nu$. This changes
          the function $h$ to $h+\nu$ and, by applying  (18), does not affect the
          truth of the result we seek. So we can suppose that $h$ is positive
          and this means that ${\cal V}$ is a positive line bundle over
          $M$. Thus we can embed $M$ as a projective variety in some
          $\bC\bP^{\mu}$, with ${\cal V}^{s}= {\cal O}(1)$, and there will
          be no loss in generality in supposing that $s=1$.
          
          Now consider our $\bC^{*}$-invariant scheme $X_{0}\subset \bC\bP^{N}$.
          Clearly we get a corresponding scheme $Z$ inside $M$, fibering
          as $\pi:Z\rightarrow \bC\bP^{1}$  with fibre $X_{0}$. For any $k$ we can identify
          $H^{0}(Z;{\cal V}^{k})$ with the sections of the vector bundle
          $\pi_{*}({\cal V}^{k})$ over $\bC\bP^{1}$. Now take the eigenspace
          decomposition
          $$  H^{0}(X_{0}, {\cal O}(k)) = \bigoplus V_{i}, $$
          where $A_{k}$ acts as $w_{i}$ on $U_{i}$. Then
          $$ \pi_{*}({\cal V}^{k})= \bigoplus U_{i} \otimes {\cal O}(w_{i}).
          $$
          It follows that $$ {\rm dim}\ H^{0}(Z,{\cal V}^{k}) = \sum {\rm dim}\
           U_{i} (w_{i}+1)
          = {\rm Tr} A_{k} + {\rm dim }H^{0} (X_{0}, {\cal O}(k)).$$
          We apply (18) to the projective scheme $Z\subset M \subset \bC\bP^{M}$.
          This shows that
          $${\rm Tr} A_{k} \sim D k^{n+1} $$
          where $D$ is the degree of $Z$. This degree is given by integrating
          $\Omega^{n+1}$ over the cycle $\vert Z\vert$, for any smooth
          form $\Omega$ on $M$ representing $c_{1}({\cal V})$. Taking the
          form given by our construction above, taking $S$ to be the
          smooth points in $\vert Z\vert$ and taking due account of multiplicity,
          we see that the degree is given by the integral $I$, and our
          result follows.

          Note that the same argument (taking $B$ to be $\bC\bP^{r}$) shows
          that for any positive integer $r$
          \begin{equation} \label{eq:genint}   {\rm Tr} A_{k}^{r} \sim k^{n+r}
          \int_{\vert X_{0}\vert} h^{r} d\mu_{FS}\ \ ,\ \ 
          {\rm Tr} \uA_{k}^{r} \sim  k^{n+r}\int_{\vert X_{0}\vert} (h-\hat{h})^{r}
          d\mu_{FS}. \end{equation}
          Thus the invariant $N_{p}$ is the $L^{p}$ norm of $h-\hat{h}$
          on $\vert X_{0}\vert$. One consequence of this is that we can
          extend Theorem 2 to the case $q=1, p=\infty$, defining $N_{\infty}({\cal
          X})= \Vert h -\hat{h}\Vert_{L^{\infty}(\vert X_{0}\vert)}$.
          It would be interesting to extend Theorem 2 to general real exponents
          $p,q$ with
          $$  N_{p} = \Vert h-\hat{h}\Vert_{L^{p}(\vert X_{0}\vert)}. $$

\end{document}